\documentclass[11pt,
]{article}                          

\usepackage{amssymb,amsfonts,amstext,amsmath,amsthm,latexsym,mathrsfs}
\usepackage{marginnote}
\usepackage{mathbbol}

\usepackage{mparhack}
\makeatletter
\if@twoside%
\else\fi%
\makeatother  
\input{cog.sty}

\begin{document}

\title
{Countable OD sets of reals belong to the ground model%.
}

\author 
{
Vladimir~Kanovei\thanks{IITP RAS and MIIT,
  Moscow, Russia, \ {\tt kanovei@googlemail.com} --- contact author. 
}  
\and
Vassily~Lyubetsky\thanks{IITP RAS,
  Moscow, Russia, \ {\tt lyubetsk@iitp.ru} 
}
}

\date 
{\today}

\maketitle

\begin{abstract}
It is true in the Cohen, random, dominaning, and 
Sacks generic extensions, that every countable 
ordinal-definable set of reals belongs to the
ground universe.
Stronger results hold in the Solovay model.
\end{abstract}

\parf{Introduction}

It is known from descriptive set theory that countable 
definable sets of reals have properties inavailable for 
arbitrary sets of reals of the same level of definability. 
Thus all elements of a countable $\id11$ set of 
reals are $\id11$ themselves while an uncountable 
$\id11$ set does not necessarily contan a $\id11$ 
real. 
This difference vanishes to some extent at higher levels 
of projective hierarchy, as it is demonstrated that 
some non-homogeneous forcing notions lead 
to models of \zfc\ 
with countable $\ip12$ non-empty sets of reals with 
no $\od$ (ordinal-definable) elements \cite{kl:cds}\snos
{The model presented in \cite{kl:cds} was obtained via 
the countable product of Jensen's minimal $\id13$ real 
forcing \cite{jenmin}. 
Such a product-forcing model was earlier considered 
by Enayat~\cite{ena}.}, 
and such a set can even have the form of   
a $\ip12$ \dd\Eo equivalence class \cite{kl:dec}. 

On the other hand, one may expect that homogeneous forcing 
notions generally yield opposite results. 
We prove the following theorems.

\bte
\lam{mt}
Let\/ $a$ be one of the following generic reals over 
the universe\/ $\rV{:}$ 
\ben
\Renv
\itlb{mt1}
a Cohen-generic real over\/ $\rV\,;$

\itlb{mt2}
a Solovay-random real over\/ $\rV\,;$

\itlb{mt3}
a dominating-forcing real over\/ $\rV\,;$

\itlb{mt4}
a Sacks (perfect-set generic) real over\/ $\rV\,.$
\een
Then it is true in\/ $\rV[a]$ that if\/ $X\sq\dn$ is 
a countable\/ $\od$ set then\/ $X\in\rV$.
\ete

\bte
\lam{mts}
{\rm(i)}
It is true in the first Solovay model\snos 
{\label{fsm1}%
See Definition \ref{sm?} below on the Solovay models.
See \cite{ksol,kl,fm16} and Stern~\cite{stern} on different 
aspects of definability in the Solovay models.} 
that every non-empty \od\ countable or finite set\/ 
$\cX$
{\bf of sets of reals}
contains an\/ \od\ element, and hence 
consists of\/ $\od$ elements  
as the notion of being \od\ is \od\ itself.

\addtocounter{footnote}{-1}

{\rm(ii)}
It is true in the second Solovay model\snom\  
that every non-empty \od\ countable or finite set\/ 
$\cX$
of any kind,
%{\bf of sets of reals}
contains an\/ \od\ element, and hence 
consists of\/ $\od$ elements, by the same reason. 
\ete

Regarding (ii), Theorem 4.8
in Caicedo and Ketchersid \cite{cai}
contains a similar result under a different \dd\ac incompatible
hypothesis on the top of $\zf+\dc$. 

One may expect such theorems to be true in any 
suitably homogeneous generic models. 
However it does not seem to be an easy task 
to manufacture a proof of sufficient
degree of generality, 
because of various {\sl ad hoc\/} arguments lacking 
a common denominator, which 
we have to make use of,
specifically for the Cohen, random, and 
dominating cases of Theorem~\ref{mt}, and 
a totally different argument used for Theorem~\ref{mts}.  

To explain the method of the proof of Theorem~\ref{mt} 
in parts \ref{mt1},  \ref{mt2},  \ref{mt3} 
(the Sacks case is quite elementary), 
let $T$ be a name of a counterexample. 
We pick a pair of reals $a,b$, each being generic over the 
ground set universe $\rV$, and satisfying $\rV[a]=\rV[b]$. 
Then the interpretations $\bint Ta\,,\:\bint Tb$ 
of $T$ resp.\ via $a$ and via $b$ coincide 
as each of them is defined by the same formula 
(with ordinals) in the same universe: $\bint Ta=\bint Tb$. 
In the same time, the pair $\ang{a,b}$ is a 
product generic pair over a suitable 
\rit{countable} model $\gM$, 
or close to be such in the sense that at least 
$\gM[a]\cap\gM[b]\cap\dn\sq\gM$. 
However $\bint Ta\sq\gM[a]$ and $\bint Tb\sq\gM[b]$, so 
in fact $\bint Ta=\bint Tb\sq\gM$, as required.

This scheme works rather transparently in the Cohen 
(Section~\ref{coh}) 
and Solovay-random (Section~\ref{ram}) 
cases, but contains a couple of 
nontrivial lemmas 
(\ref{cH1} and especially \ref{cH2} with a 
lengthy proof) in the dominating case 
(Section~\ref{domi}).  

We add an alternative and rather elementary proof 
for the Cohen and Solovay-random cases 
(Section~\ref{c+r}), which makes use of some old 
folklore results related to degrees of reals 
in those extensions over the ground model.  
We finish in Section \ref{solm} with a proof 
of Theorem~\ref{mts}.

\parf{Cohen-generic case}
\las{coh}

Here we prove Case \ref{mt1} of Theorem~\ref{mt}.
We begin with some notation and a couple of 
preliminary lemmas.

Assume that $u,v\in\dn\cup\bse$ are dyadic sequences,
possibly of different 
(finite or infinite) length. 
We let $u\oq v$ (the termwise action of $u$ on $v$) 
be a dyadic sequence defined so that 
$\dom{u\oq v}=\dom v$ 
(independently of the length $\dom u$ of $u$)
and if $j<\dom v$ then 
$$
(u\oq v)(j)=
\left\{
\bay{rcl}
1-v(j)&, \text{ whenever}& j< \dom{u}\land u(j)=1\:,\\[1ex]
v(j)&, \text{ otherwise}& . 
\eay
\right.
$$
In particular, if $z\in\dn\cup\bse$ then 
$x\mto z\oq x$ ($x\in\dn$) is a homeomorphism of $\dn$ while 
$p\mto z\oq p$ ($p\in\bse$) is an order automorphism of $\bse$.

Let $\koh=\bse$ be the Cohen forcing.

\ble
\lam{aL}
Let\/ $\gM$ be a transitive model of a large 
fragment of\/ $\zfc$.
Then
\ben
\renu
\itlb{aL0}
if a pair\/ $\ang{a,b}\in\dn\ti\dn$ is\/
\dd{(\koh\ti\koh)}generic over\/ $\gM$ then\/
$\gM[a]\cap\gM[b]=\gM$ --- 
{\rm this is a well-known theorem on product forcing}$;$

\itlb{aL1}
if a pair\/ $\ang{a,b}\in\dn\ti\dn$ is\/
\dd{(\koh\ti\koh)}generic over\/ $\gM$ then
so is the pair\/ $\ang{a,a\aq b}\;;$

\itlb{aL2}
if\/ $\gM$ is countable and\/ $p,q\in\koh$ then
there are reals\/ $a,b\in\dn$, \dd\koh generic over\/ 
$\rV$ and such that\/ 
$p\su a$, $q\su b$, 
$\rV[a]=\rV[b]$,
%, $x\aq y\in\rV$,
and the pair\/ $\ang{a,b}$ 
is\/ \dd{(\koh\ti\koh)}generic over\/ $\gM$.
\een
\ele
\bpf
\ref{aL1}
Otherwise there is a condition $\ang{p,q}\in \koh\ti\koh$
with $\dom p=\dom q$, 
which forces the opposite over $\gM$.
By the countability, there is a real $a\in\dn$ in $\rV$ 
\dd\koh generic over $\gM$, with $p\su a$; 
$\gM[a]$ is a set in $\rV$.
Let $r=p\aq q$ and let $c\in\gM$ be
\dd\koh generic over $\gM[a]$, with $r\su c$.
Then $b=a\aq c$ is \dd\koh generic over $\gM[a]$
by obvious reasons, $c=a\aq b$,
and $q=p\aq r\su b=a\aq c$.
Finally $\ang{a,b}$ is 
\dd{(\koh\ti\koh)}generic over $\gM$  
by the product forcing theorem, a contradiction.

\ref{aL2}
Assuming wlog that $\dom p=\dom q$, we let $r=p\aq q$. 
Once again, there is a real $c\in\dn$ in $\rV$, 
\dd\koh generic over $\gM$, with $r\su c$.
Let $a\in\dn$ be \dd\koh generic over $\rV$, hence
over $\gM[c]$, too, and satisfying $p\su a$. 
Then the real $b=c\aq a$ 
is \dd\koh generic over $\rV$ (since $c\in\rV$),   
$\rV[b]=\rV[a]$, and $q=r\aq p\su b$. 
 
Finally the pair $\ang{a,c}$ is 
\dd{(\koh\ti\koh)}generic over $\gM$  
by the product forcing theorem, therefore
$\ang{a,b}=\ang{a,a\aq c}$ is 
\dd{(\koh\ti\koh)}generic over $\gM$ by \ref{aL1}.
\epf

\vyk{
Let $\doa$ be the canonical \dd\koh name for the 
\dd\koh generic real.

Let  
$\dal,\dar$ be canonical \dd{(\koh\ti\koh)}names for the left, 
resp., right of the terms of a\/ \dd{(\koh\ti\koh)}generic pair 
of reals\/ $\ang{\alev,\apra}\in\dn\ti\dn$.

Let $\uv$ be a name for the ground model 
(of ``old'' sets).
}

\bpf[Theorem~\ref{mt}, case \ref{mt1}]
Let $a_0\in\dn$ be a real
\dd\koh generic over the universe $\rV$. 
First of all, note this: 
it suffices to prove that 
(it is true in\/ $\rV[a_0]$ that) if\/ $Z\sq\dn$ is 
a countable\/ $\od$ set then\/ $Z\sq\rV$.
Indeed, as the Cohen forcing is homogeneous, any 
statement about sets in $\rV$, the ground model, 
is decided by the weakest condition.

Thus let $Z\sq\dn$ be a countable\/ $\od$ set in $\rV[a_0]$.

{\ubf Suppose to the contrary that $Z\not\sq\rV$}.

There is a formula $\vpi(z)$ 
with an unspecified ordinal $\ga_0$ as a parameter, 
such that $Z=\ens{z\in\dn}{\vpi(z)}$ in $\rV[a_0]$, and then 
there is a condition $p_0\in\koh$ such that $p_0\su a_0$ and 
$p_0$ \dd\koh forces that 
$\ens{z\in\dn}{\vpi(z)}$ is a countable set and 
(by the contrary assumption) 
also forces $\sus z\:(z\nin\rV\land\vpi(z))$.

There is a sequence $\sis{t_n}{n<\om}\in\rV$ of \dd\koh names, 
such that if $x\in\dn$ is Cohen generic and $p_0\su x$ then 
it is true in $\rV[x]$ that 
$\ens{z\in\dn}{\vpi(z)}=\ens{\bint{t_n}x}{n<\om}$, 
where $\bint{t}x$ is the interpretation of a \dd\koh name $t$ 
by a real $x\in\dn.$
Let $T\in\rV$ be the canonical \dd\koh name for 
$\ens{\bint{t_n}{\dot a}}{n<\om}$.
Thus we assume that
\benq
\nenu
\itlb{enu1}
$p_0$ \dd\koh forces, over $\rV$, that 
$\bint{T}{\dot a}=
%\ens{\bint{t_n}{\dot a}}{n<\om}=
\ens{x\in\dn}{\vpi(x)}
\not\sq\uv$, 
%\eqno(1)
\een
where $\doa$ is the canonical \dd\koh name for the 
\dd\koh generic real,
and $\uv$ is a name for the ground model (of ``old'' sets).
%
%where $\La$, the empty string, is the weakest condition in $\koh$.

We continue towards 
{\ubf getting a contradiction from \ref{enu1}}. 
Pick a regular cardinal $\ka>\al_0$, sufficiently large for 
the set $\hk\ka$ to contain $\ga_0$ and all names $t_n$ and $T$.
%and  for Lemma~\ref{aL}\ref{aL0},\ref{aL1} to be true for
%$\gM=\hk\ka$. 
%
Consider a countable elementary submodel $\gM$ of $\hk\ka$ 
containing $\ga_0$, all $t_n$, $T$. 
Let $\pi:\gM\to\gM'$ be the Mostowski collapse onto a 
transitive set $\gM'$.
As $\koh$ is countable, we have
$\pi(\koh)=\koh$, $\pi(t_n)=t_n$, $\pi(T)=T$, so $T\in\gM'$.

Now pick reals $a,b\in\dn$ \dd\koh generic over $\rV$
by Lemma~\ref{aL}\ref{aL2},
%(with respect to $\gM'$),
such that $p_0\su a$, $p_0\su b$, $\rV[a]=\rV[b]$,
and the pair
$\ang{a,b}$ is \dd{(\koh\ti\koh)}generic over $\gM'$.
In particular, as $\rV[a]=\rV[b]$, 
% and $p_0\su x$, $p_0\su y$, 
we have $\bint{T}{a}=\bint{T}{b}\not\sq\rV$ by 
\ref{enu1}.
On the other hand,
%Note that
$\gM'[a]\cap\gM'[b]\sq\gM'$
by Lemma~\ref{aL}\ref{aL0}, therefore
$\bint{T}{a}\cap\bint{T}{b}\sq\gM'[a]\cap\gM'[b]\sq\gM'\sq\rV$,
contrary to the above.\vom

\epF{Theorem~\ref{mt}, case \ref{mt1}}

\parf{Solovay-random case}
\las{ram}

Here we prove Case \ref{mt2} of Theorem~\ref{mt}.

Let $\jla$ be the standard probability Lebesgue measure 
on $\dn.$
The Solovay-random forcing $\srf$ consists of all 
trees $\tau\sq\bse$ with  
no endpoints and no isolated branches, and such that the set  
$[\tau]=\ens{x\in\dn}{\kaz n\,(x\res n\in\tau)}$ has 
positive measure $\jla([\tau])>0$.
%Unlike $\koh$, t
The forcing $\srf$ 
depends on the ground model, so that ``random over 
a model $\gM$'' will mean 
``\dd{(\srf\cap\gM)}generic over $\gM$''. 

\ble
[trivial in the Cohen case]
\lam{bJ}
If\/ $\gM\sq\gN$ are TM of a large 
fragment of\/ $\zfc$, and\/ $a\in\dn$ is random over\/ 
$\gN$ then\/ $a$ is random over\/ $\gM$, too.
\ele
\bpf
It suffices to prove that if $A\in\gM$ is a maximal 
antichain in $\srf\cap\gM$ then $A$ remains such in 
$\srf\cap\gN$, which is rather clear since being 
a maximal antichain in $\srf$ amounts to 
1) countability, 
2) pairwise intersections being null sets
(those of \dd\jla measure $0$), 
and 3) the union being a co-null set.
\epf

Unlike the Cohen-generic case, a random pair of reals is
{\ubf not} a \dd{(\srf\ti\srf)}generic pair.
The notion of a random pair 
is rather related to forcing by closed sets in
$\dn\ti\dn$ (or trees which generate them, or equivalently
Borel sets)
of positive product measure (non-null).
This will lead to certain changes of arguments, with respect 
to the Cohen-generic case of Section~\ref{coh}.

We'll make use of the following known characterisation of
random pairs.

\bpro
\lam{bK}
Let\/ $\gM$ be a transitive model of a large 
fragment of\/ $\zfc$, and\/ $a,b\in\dn$.
Then the following 
three assertions are equivalent$:$\vom

$1)$ 
the pair\/ $\ang{a,b}$ is a random pair over\/ $\gM\;;$\vom

$2)$ 
$a$ is random over\/ $\gM$ and\/
$b$ is random over\/ $\gM[a]\;;$\vom

$3)$ 
$b$ is random over\/ $\gM$ and\/
$a$ is random over\/ $\gM[b]\;.$\qed
\epro

\ble
\lam{bL}
Let\/ $\gM$ be a transitive model of a large 
fragment of\/ $\zfc$.
Then
\ben
\renu
\itlb{bL0}
if a pair\/ $\ang{a,b}\in\dn\ti\dn$ is\/
random over\/ $\gM$ then\/
$\gM[a]\cap\gM[b]\cap\dn\sq\gM\,;$

\itlb{bL1}
if a pair\/ $\ang{a,b}\in\dn\ti\dn$ is\/
random over\/ $\gM$ then
so is the pair\/ $\ang{a,a\aq b}\,;$

\itlb{bL2}
if\/ $\gM$ is countable and\/ $\tau\in\srf$ then
there are reals\/ $a,b\in[\tau]$, random over\/ 
$\rV$, such that\/ $\rV[a]=\rV[b]$,
%, $x\aq y\in\rV$,
and the pair\/ $\ang{a,b}$ is\/ random over\/ $\gM$.
\een
\ele
\bpf
\ref{bL0}
This is somewhat more difficult than in the
Cohen-generic case of Lemma~\ref{aL}\ref{aL0}.
Assume towards the contrary that
$x\in\gM[a]\cap\gM[b]\cap\dn$ but $x\nin M$.
The random forcing admits continuous reading of real
names, meaning that there are continuous maps
$f,g:\dn\to\dn$, coded in $\gM$ and such that $x=f(a)=g(b)$.
Let the contrary assumption be forced by a Borel set
$P\sq\dn\ti\dn$ of positive product measure, 
coded in $\gM$ and containing $\ang{a,b}$; 
in particular, $P$ 
(random pair)-forces that $f(\dal)=g(\dar)$.\snos
{$\dal,\dar$ are canonical names for the left, 
resp., right of the terms of a random pair.}
By the Lebesgue density theorem, we can wlog assume that 
every point $\ang{x,y}\in P$ has density $1$. 

We claim that $f(x)=g(y)$ for all $\ang{x,y}\in P$. 
Indeed if $\ang{x_0,y_0}\in P$ and $f(x_0)\ne g(y_0)$ then 
say $f(x_0)(n)=0\ne g(y_0)(n)=1$ for some $n$. 
As $f,g$ are continuous, there is a nbhd $Q$ of 
$\ang{x_0,y_0}$ in $P$ such that $f(x)(n)=0\ne g(y)(n)=1$
for all $\ang{x,y}\in Q$. 
But $Q'$ is a non-null set by the density 1 
assumption. 
It follows that $Q$ forces that $f(\dal)\ne g(\dar)$, a
contradiction.

Let a \rit{cell} be any Borel set $Q\sq P$ such that $f,g$ 
are constant on $Q$, that is, there is a real $r$ such that 
$f(x)=g(y)=r$ for all $\ang{x,y}\in Q$.
Note that in this case, if $Q$ is non-null then $Y$ forces 
$f(\dal)=g(\dar)=r\in\gM$, therefore to prove \ref{bL0} it 
suffices to show the existence of a non-null cell $Q\sq P$. 

Let $P_x=\ens{y}{\ang{x,y}\in P}$ and 
$P^y=\ens{x}{\ang{x,y}\in P}$, cross-sections. 
By Fubini, the sets $X=\ens{x}{\jla(P_x)>0}$  
and $Y=\ens{y}{\jla{(P^y\cap X)}>0}$ are non-null. 
Let $y_0\in Y$ and let $X'= P^{y_0}\cap X$, a non-null set. 
By construction, if $x\in X'$ then the cross-section $P_x$ 
is non-null, and hence 
$Q=\ens{\ang{x,y}\in P}{x\in X'}$ is non-null by Fubini. 
We claim that $Q$ is a cell. 
Indeed suppose that $\ang{x,y}\in Q$. 
Then $x\in X'$, therefore  
$\ang{x,y_0}\in P$, and 
we have $f(x)=g(y_0)$ by the above claim. 
However $\ang{x,y}\in P$, 
hence similarly $g(y)=f(x)$. 
Thus $g(y)=f(x)=g(y_0)=\text{Const}$ on $Q$, as required.

\ref{bL1}
The contrary assumption implies the existence (in $\gM$)
of a \rit{non-null} Borel set $P\sq\dn\ti\dn$ and 
a \rit{null} Borel set $Q\sq\dn\ti\dn$ such that the 
map $\ang{x,y}\mto\ang{x,x\aq y}$ maps $P$ into $Q$. 
However this map is obviously measure-preserving, 
a contradiction.

\ref{bL2}
The set $P=\ens{\ang{x,x\aq y}}{x,y\in[\tau]}$ is non-null,
hence,  by Fubini, so is the projection $Y=\ens{y}{\jla(P^y)>0}$, 
where $P^y=\ens{x}{\ang{x,y}\in P}$, as above. 
Let, in $\rV$, $y\in Y$ be random over $\gM$.
Then $P^y$ is non-null, so we can pick a real $a\in P^y$
random over $\rV$ hence, over $\gM[y]$, too.
Then the pair $\ang{a,y}$ belongs to $P$ and is random over
$\gM$ by Proposition~\ref{bK}.
Let $b=a\aq y$.
It follows by \ref{bL1} that the pair $\ang{a,b}$
is random over $\gM$ as well.
And $a,b\in[\tau]$ by construction. 
Finally $b$ is random over $\rV$ since so is $a$ while
$y\in\rV$. 
\epf

%\newpage

\bpf[Theorem~\ref{mt}, case \ref{mt2}]
As above (the Cohen case), 
the {\ubf contrary assumption} 
leads to a formula $\vpi(z)$ with $\ga_0\in\Ord$ as a parameter, 
a condition $\tau_0\in\srf$ in $\rV$
which \dd\srf forces, over $\rV$, that the set 
$\ens{z\in\dn}{\vpi(z)}$ is countable and 
%(by the contrary assumption) 
$\sus z\:(z\nin\uv\land\vpi(z))$,
a sequence $\sis{t_n}{n<\om}\in\rV$ of
\dd\srf names for reals in $\zn$, 
and a canonical \dd\srf name $T\in\rV$ for 
$\ens{\bint{t_n}{\dot a}}{n<\om}$, 
such that 
\vyk{
Let $a_0\in\dn$ be a random real over the universe $\rV$. 
As above, 
it suffices to prove that if 
$\vpi(z)$ is a formula with $\ga_0\in\Ord$ as a parameter 
and $Z=\ens{z\in\dn}{\vpi(z)}$ is 
a countable $\od$ set in $\rV[a_0]$ then $Z\sq\rV$.

{\ubf Suppose to the contrary that this fails.}
Then there is a condition $\tau_0\in\srf$ in $\rV$
such that $a_0\in[\tau_0]$ and 
$\tau_0$ \dd\srf forces, over $\rV$, that 
$\ens{z\in\dn}{\vpi(z)}$ is a countable set and 
(by the contrary assumption) 
also forces $\sus z\:(z\nin\uv\land\vpi(z))$.
%where $\uv$ is a name for the ground model.

There is a sequence $\sis{t_n}{n<\om}\in\rV$ of
\dd\srf names for reals, 
such that if $x\in[\tau_0]$ is random then 
it is true in $\rV[x]$ that 
$\ens{z\in\dn}{\vpi(z)}=\ens{\bint{t_n}x}{n<\om}$, 
where $\bint{t}x$ is the interpretation of a \dd\srf name $t$ 
by a real $x\in\dn.$
Let $T\in\rV$ be the canonical \dd\srf name for 
$\ens{\bint{t_n}{\dot a}}{n<\om}$.
Thus we assume that
}
\vyk{
\benq
\nenu
%\atc
\itlb{enu2}
$\tau_0$ \dd\srf forces, over $\rV$, that 
$\bint{T}{\dot a}=
%\ens{\bint{t_n}{\dot a}}{n<\om}=
\ens{x\in\dn}{\vpi(x)}
\not\sq\uv$, 
%\eqno(1)
\een
where $\doa$ is the canonical \dd\srf name for the 
\dd\srf generic real.
}

\benq
\nenu
\itlb{enu2}
if $x\in[\tau_0]$ is a random real over $\rV$, 
%compatible with ${p_0}$ 
then it is true in $\rV[x]$ that 
$$
\ens{z\in\dn}{\vpi(z)}=\ens{\bint{t_n}x}{n<\om}=\bint Tx 
\not\sq\rV\,.
$$
\een

%
%where $\La$, the empty string, is the weakest condition in $\koh$.

Pick a regular cardinal $\ka>\al_0$, sufficiently large for 
the set $\hk\ka$ to contain $\ga_0$ and all names $t_n$ and $T$.
%and  for Lemma~\ref{aL}\ref{aL0},\ref{aL1} to be true for
%$\gM=\hk\ka$. 
%
Consider a countable elementary submodel $\gM$ of $\hk\ka$ 
containing $\ga_0$, all names $t_n$ and $T$, and $\srf$. 
Let $\pi:\gM\to\gM'$ be the Mostowski collapse onto a 
transitive set $\gM'$.
Unlike the Cohen case, the set $\srf'=\pi(\srf)$ 
is equal to $\srf\cap\gM'$, just the random forcing in 
$\gM'$, but still $\pi(t_n)=t_n$ for all $n$, 
since by the ccc property of $\srf$ we can assume that 
$t_n$ is a hereditarily countable set, 
and accordingly $\pi(T)=T$. 

Pick reals $a,b\in[\tau_0]$ random over $\rV$
by Lemma~\ref{bL}\ref{bL2},
%(with respect to $\gM'$),
such that $\rV[a]=\rV[b]$, and the pair
$\ang{a,b}$ is random over $\gM'$.
As $\rV[a]=\rV[b]$, 
% and $p_0\su x$, $p_0\su y$, 
we have $\bint{T}{a}=\bint{T}{b}\not\sq\rV$ by 
\ref{enu2}.
But $\gM'[a]\cap\gM'[b]\sq\gM'$
by Lemma~\ref{aL}\ref{aL0}, therefore
$\bint{T}{a}\cap\bint{T}{b}\sq\gM'[a]\cap\gM'[b]\sq\gM'\sq\rV$, 
and we get a contratiction required.\vom 

\epF{Theorem~\ref{mt}, case \ref{mt2}}

\parf{Cohen and random cases: a different proof}
\las{c+r}

Here we present a shorter proof of Cases \ref{mt1} 
and \ref{mt2} of Theorem~\ref{mt}, based on the 
following lemmas.

\ble
\lam{addc}
Let\/ $a\in\dn$ be Cohen-generic over a transitive 
model\/ $\gM$, and\/ $b\in\dn\cap\gM[a]$, a real in 
the extension. 
Then
\ben
\renu
\itlb{addc1}
either\/ $b\in\gM$ or there is a real\/ $b'\in\dn$, 
Cohen-generic over\/ $\gM$ and satisfying\/ 
$\gM[b]=\gM[b']\;;$

\itlb{addc2}
either\/ $\gM[b]=\gM[a]$ or\/ $\gM[a]$ is a 
Cohen-generic extension of\/ $\gM[b]$.\qed
\een
\ele

\ble
\lam{addr}
Let\/ $a\in\dn$ be random over a transitive 
model\/ $\gM$, and\/ $b\in\dn\cap\gM[a]$, a real in 
the extension. 
Then
\ben
\renu
\itlb{addr1}
either\/ $b\in\gM$ or there is a real\/ $b'\in\dn$, 
random over\/ $\gM$ and satisfying\/ 
$\gM[b]=\gM[b']\;;$

\itlb{addr2}
either\/ $\gM[b]=\gM[a]$ or\/ $\gM[a]$ is a 
random extension of\/ $\gM[b]$.\qed
\een
\ele

The lemmas are known in set theoretic folklore, yet 
we are not able to suggest any reference. 
In particular Lemma~\ref{addc}\ref{addc2} is 
rather simple  
on the base on general results on intermediate models
by Grigorieff~\cite{gri} since any subforcing 
of the Cohen forcing 
either is trivial or is equivalent to Cohen forcing. 

\bpf
[Theorem~\ref{mt}, case \ref{mt1}, from Lemma~\ref{addc}]
In $\gM[a]$, let $b$ belong to a countable \od\ set 
$X=\ens{x\in\dn}{\vpi(x)}$, where $\vpi$ is a formula   
containing ordinals.
As $b\in\gM[a]$, there is a Borel function $f$, 
coded in $\gM$, such that $b=f(a)$. 
We have to prove that $b\in\gM$.
Let 
%$\dob$ be any \dd\koh name for $b$ and 
$\doa$ be a canonical \dd\koh name for the generic real.

We have two cases, by Lemma~\ref{addc}\ref{addc2}.\vom 

{\it Case 1\/}: $\gM[b]=\gM[a]$. 
Then there is a Borel function $g$, coded in $\gM$, 
such that $a=g(b)$. 
There is a Cohen condition $u\in\koh$ which 
satisfies $u\su a$ and forces 
$\doa=g(f(\doa))$, $\vpi(f(\doa))$, and the sentence 
``$\ens{x\in\dn}{\vpi(x)}$ is countable''.

Now, the set $A$ of all reals $a'\in\dn$, 
Cohen-generic over $\gM$ and satisfying 
$u\su a'$ and $\gM[a']=\gM[a]$, 
belongs to $\gM$ and definitely is 
uncountable in $\gM$. 
If $a'\in A$ then $f(a')$ satisfies $\vpi(f(a'))$ in 
$\gM[a']=\gM[a]$ and hence belongs to $X$. 
Furthermore if $a'\ne a''\in A$ then $f(a')\ne f(a'')$ 
since $a'=g(f(a')$ and $a''=g(f(a'')$. 
We conclude that $X$ is uncountable, a contradiction.\vom

{\it Case 2\/}: $\gM[a]$ is a Cohen-generic extension 
of $\gM[b]$. 
Let $\psi(x)$ be the formula saying: 
``$x\in\dn$ and $\koh$ forces $\vpi(\dotx)$, where $\dotx$ 
is a canonical \dd\koh name for $x$ in any transitive 
graund model containing $x$. 
As $\koh$ is a homogeneous forcing notion, the set 
$Y=X\cap\gM[b]$ coincides with the set 
$\ens{x\in\dn}{\psi(x)}$ defined in $\gM[b]$, and 
$b\in Y$. 
Finally $\gM[b]$ is a Cohen extension of $\gM$ 
by Lemma~\ref{addc}\ref{addc1} 
(or else just $b\in\gM$), 
and it remails to apply the result in Case 1 to $Y$. 
\epf 

\bpf
[Theorem~\ref{mt}, case \ref{mt2}, from Lemma~\ref{addr}]
Similar. 
\epf 

It is really temptating to prove the dominating case 
of the theorem by this same rather simple method. 
However we cannot establish any result similar to 
lemmas \ref{addc}, \ref{addr} for dominating forcing. 
Some relevant results by Palumbo \cite{paldis,pal} 
fall short of what would be useful here. 
Generally, a remark in \cite[Section 4]{pal} casts 
doubts that even claims \ref{addr1} of the lemmas 
hold for dominating-generic extensions in any useful 
form. 
This is why we have to process the dominating case 
of Theorem~\ref{mt} the hard way in the next section.

\parf{Dominating case}
\las{domi}

Here we prove Case \ref{mt3} of Theorem~\ref{mt}.

Let $\dZ=\ans{\dots,-2,-1,0,1,2,\dots}$,                   
integers of both signs.

We let the dominating forcing $\df$
consist of all pairs $\ang{n,f}$ such that $f\in\dZ^\om$
(that is, $f$ is an infinite sequence of integers) and 
$n<\om$.
We order $\df$ so that $\ang{n,f}\le\ang{n',f'}$ 
(the bigger is stronger) iff
$n\le n'$, $f\res n=f'\res n$, and $f\le f'$ componentwise,
that is, $f(k)\le f'(k)$ holds for all $k<\om$.\snos
{This slightly differs from the standard definition, as \eg\
in Bartoszy\'nski -- Judah \cite[3.1]{bj} where $f\in\bn$.
The difference does not change any forcing properties,
but leads to a more friendly setup since
$\df$ as defined here 
is a group under componentwise addition.}

A modified version $\df'$ consists of all pairs $\ang{u,h}$, 
where $u\in\zse$, $h\in\zn$. 
Each such pair is 
identified with the pair 
$\ang{\dom u,u\we h}\in\df$, where $\we$  denotes the 
concatenation, and the order on $\df'$ is induced by this 
identification.

\bdf
\lam{ddf}
If $G\sq\df$ is a generic filter then
$a_G=\bigcup_{\ang{n,f}\in G}f\res n$ belongs to $\dZ^\om$;
we call $a_G$ {\it a dominating-generic real\/}.
More exactly, if $\gM$ is a transitive model and a set
$G\sq\df\cap\gM$ is \dd{(\df\cap\gM)}generic over $\gM$
then say that $a_G$ is
a {\it dominating-generic\/} (\dg, in brief)
real over\/ $\gM$. 
\edf

\bre
Unfortunately there is no result similar to 
Proposition~\ref{bK} for the dominating forcing, 
since if $a$ is a \dg\ real over $\gM$ and 
$b$ is a \dg\ real over $\gM[a]$ then $a$ is definitely 
not \dg\ over $\gM[b]$. 
This will make our arguments here somewhat more complex 
than in the Solovay-random section.
\ere

If $u,v$ are finite or infinite sequences of integers in $\dZ$ 
then let $u\ar v$ be a sequence defined by componentwise sum,
so that $\dom{(u\ar v)}=\dom v$ 
(independently of the length $\dom u$)
and if $j<\dom v$ then $(u\ar v)(j)=u(j)+v(j)$.
If in addition $\dom u=\dom v$ then $u\am v$ is defined
similarly.

For instance $f\ar g$ and $f\am g$ are defined for all $f,g\in\zn.$ 

\vyk{
In addition, if 
$\ang{n,f}\in\df$ then let $u\ar \ang{n,f}=\ang{n,u\ar f}$.
The map $\ang{n,f}\mto u\ar \ang{n,f}$ is an order
automorphism of $\df$ provided $u\in\zn\cup\zse$.
}

\ble
%[trivial in Cohen-generic case]
\lam{cI}
If\/ $\gM\sq\gN$ are TM of a large 
fragment of\/ $\zfc$, and\/ $a\in\zn$ is \dg\ 
over\/ $\gN$ then\/ $a$ is \dg\ over\/ $\gM$, too.
\ele
\bpf
It suffices to prove that if $A\in\gM$ is a maximal 
antichain in $\df\cap\gM$ then $A$ remains such in 
$\df\cap\gN$.
Note that $A$ is countable in $\gM$
since $\df$ is a ccc forcing, therefore $A$ is effectively
coded by a real $r\in\gM$ so that being a maximal antichain
is a $\ip11$ property of $r$. 
It remains to refer to the Mostowski absoluteness theorem.
\epf

\ble
\lam{cJ}
If\/ $\gM$ is a TM of a large 
fragment of\/ $\zfc$, $h\in\gM\cap\zn,$
and\/ $a\in\zn$ is a \dg\ real
over\/ $\gM$ then\/ $a\ar h$, $a\am h$ are
\dg\ over\/ $\gM$, too.
\ele
\bpf
The maps $\ang{n,f}\mto\ang{n,f\ar h}$ and
$\ang{n,f}\mto\ang{n,f\am h}$ are order-automorphisms of
$\df\cap\gM$ in $\gM$.
\epf

\ble
%[trivial in Cohen-generic case]
\lam{cH1}
If\/ $\gM$ is a TM of a large fragment of\/ $\zfc$,
$a\in\zn$ is a \dg\ real over\/ $\gM$, and\/
$b\in\zn$ is a \dg\ real over\/ $\gM[a]$,  
then\/ $\gM[a]\cap\gM[b]\cap\dn\sq\gM$.
\ele
\bpf
Otherwise the opposite is forced over $\gM[a]$ by
a condition $\ang{n,f}\in\df\cap\gM[a]$;
thus $f\in\zn\cap\gM[a]$.
To be more precise, $\ang{n,f}$
\dd{(\df\cap\gM)}forces
$\ugm[\dob]\cap\ugm[a]\cap\dn\not\sq\ugm$ 
over $\gM[a]$, where 
$\ugm$ is a suitable name for $\gM$ as a class in $\gM[a]$,
and $\dob$ is a canonical name  
for the \dg\ real over $\gM[a]$.

We claim that any other condition 
$\ang{n',f'}\in\df\cap\gM[a]$ forces the same. 
Suppose to the contrary that in fact some 
$\ang{n',f'}\in\df\cap\gM[a]$ forces 
$\ugm[\dob]\cap\ugm[a]\cap\dn\sq\ugm$ 
over $\gM[a]$.
We can wlog assume that $n'=n$ and the \dd ntails 
of $f$ and $f'$ coincide: 
$f(j)=f'(j)$ for all $j\ge n$. 
Now let $b\in\zn$ be a \dg\ real over $\gM[a]$ compatible 
with $\ang{n,f}$, that is, $b\res n=f\res n$ and $f\le b$ 
componentwise. 
Let $b'\in\zn$ be defined so that 
$b'(j)=b(j)$ for all $j\ge n$, but $b\res n=f'\res n$;      
then $b'$ is a \dg\ real over $\gM[a]$ 
compatible with $\ang{n,f'}$.
Then by construction we have 
$\gM[b]\cap\gM[a]\cap\dn\not\sq\gM$ but 
$\gM[b']\cap\gM[a]\cap\dn\sq\gM$.
However obviously $\gM[b]=\gM[b']$, a 
contradiction which completes the claim.

We conclude that if $b\in\zn$ is any \dg\ real 
over $\gM[a]$ then $\gM[b]\cap\gM[a]\cap\dn\not\sq\gM$. 
As $a$ itself is generic over $\gM$, there is a 
condition $\ang{m,h}\in\df\cap\gM$ such that 
$\gM[b]\cap\gM[a]\cap\dn\not\sq\gM$ holds whenever 
$a\in\zn$ is \dg\ over $\gM$ compatible with $\ang{m,h}$ 
and $b\in\zn$ is \dg\ over $\gM[a]$.

Now let $\ka=2^{\aleph_0}$ in $\gM$, 
and let $\la=\ka^+$ be the next cardinal in $\gM$. 
Let 
$$
\dQ=\ens{\ang{m',h'}\in\df\cap\gM}{\ang{m,h}\le\ang{m',h'}}\,.
$$ 
Consider the finite-support forcing product $\dQ^\la$ 
in $\gM$. 
A \dd{\dQ^\la}generic extension of $\gM$ has the form 
$\gN=\gM[\sis{a_\xi}{\xi<\la}]$, 
where $a_\xi\in\dn$ are pairwise \dg\ reals over 
$\gM$, compatible with $\ang{m,h}$, in particular 
$\gM[a_\xi]\cap\gM[a_\eta]=\gM$ whenever $\xi\ne\eta$. 

Consider a \dd{(\df\cap\gN)}generic extension $\gN[b]$ 
of $\gN$, so that $b\in\zn$ is a \dg\ real over $\gN$. 
Then $b$ is \dg\ over each $\gM[a_\xi]$ by Lemma~\ref{cI}.
It follows by the above that 
$\gM[b]\cap\gM[a_\xi]\cap\dn\not\sq\gM$. 
Let $z_\xi\in \gM[b]\cap\gM[a_\xi]\cap\dn\bez\gM$, 
for all $\xi<\la$. 
Note that if $\xi\ne\eta$ then $z_\xi\ne z_\eta$ 
since $\gM[a_\xi]\cap\gM[a_\eta]=\gM$, see above. 
Thus we have \dd\la many different reals in $\gM[b]$. 
However $\gM[b]$ is a CCC extension 
of $\gM$ by Lemma~\ref{cI}, and hence there cannot be more 
(in the sense of cardinality) 
reals in $\gM[b]$ than in $\gM$. 
% (\dd\ka many of them).
The contradiction ends the proof.
\epf

\ble
%[trivial in Cohen-generic case]
\lam{cH2}
If\/ $\gM$ is a TM of a large fragment of\/ $\zfc$,
$a\in\zn$ is a \dg\ real over\/ $\gM$, and\/
$b\in\zn$ is a \dg\ real over\/ $\gM[a]$,  
then\/ $\gM[b]\cap\gM[a\ar b]\cap\dn\sq\gM$.
\ele

One may want to prove the lemma by proving that 
$\ang{b,a\ar b}$ is dominating product-generic over $\gM$ 
due to the genericity of $a$. 
But in fact this is not the case. 
Indeed if $\ang{b,a\ar b}$ is dominating product-generic over 
$\gM$ then a transparent forcing argument shows that 
$a=(a\ar b)\am b$ is simply Cohen-generic over $\gM$, contrary 
to $a$ being \dg.

\bpf
%\value{equation}5
By Lemma~\ref{cJ}, $a\ar b$ is \dg\ over $\gM[a]$, 
and hence over $\gM$ by Lemma~\ref{cJ} 
Therefore the {\ubf contrary assumption} implies
a pair of \dd{(\df\cap\gM)}real names $\sg,\tau\in\gM$ such
that $\bint\sg b =\bint\tau{a\ar b}\in\dn\bez\gM$, where 
$\bint{t}b$ is the \dd binterpretation of $\sg$.

Let us present 
the two-step iterated forcing ${\dP}\in\gM$ which produces 
$\gM[a][b]$ as $\df\ast\dfp$, with $\dfp$, not $\df$, 
as the second stage. 
Then ${\dP}$ consists of all quadruples, 
or double-pairs, of the form 
$p=\ang{\ang{m_p,f_p},\ang{u_p,t_p}}=\ang{m_p,f_p,u_p,t_p}$, 
where $\ang{m_p,f_p}\in\df\cap\gM$, 
$u_p\in\zse$, and 
$t_p\in\gM$ is a \dd\df name for an element 
of $\zn$, 
%$\dZ^{\om\bez\dom{u_p}}$, 
with a suitable order. 
We shall use $\doa,\,\dob$ as canonical \dd\dP names of the 
\dg\ real over $\gM$ and \dg\ real over $\gM[a]$, 
respectively. 

By the contrary assumption, 
there is a condition $p_0=\ang{m_0,f_0,u_0,t_0}\in {\dP}$
which \dd{{\dP}}forces, over $\gM$, the formula 
$
\bint\sg\dob=\bint\tau{\doa\ar \dob}\in\dn\bez\ugm\,,
$ 
so that 
\benq
\nenu
%\setcounter{enumi}{\value{equation}}
%\stepcounter{equation}
\itlb{xx1}
if $\ang{a,b}\in\zn$ is a pair \dd\dP generic over $\gM$ 
(so $a$ is \dg\ over $\gM$ and $b$ \dg\ over $\gM[a]$)  
and compatible with $p_0$, then 
%we have 
$
\bint\sg b =\bint\tau{a\ar b}\in\dn\bez\gM\,.
$
\een 
Let $n_0=\dom{u_0}$. 
We can assume that $n_0\le m_0$; 
otherwise change $m_0$ to $n_0$. 

By simple strengthening, we find a stronger condition 
$p_1=\ang{m_1,f_1,u_1,t_1}$ in ${\dP}$, $p_1\ge p_0$, such that 
$m_0\le n_1=\dom{u_1}\le m_1$.

\bcl
\lam{xx2}
If conditions\/ $p_2=\ang{m,f,u_2,t_2}$ and 
$p_3=\ang{m,f,u_3,t_3}$ $($same\/ $m,f$$!)$ in\/ ${\dP}$ satisfy 
$p_1\le p_2$, $p_1\le p_3$, 
%and $\ang{m_1,f_1}=\ang{m_2,f_2}$, 
and in addition\/ $k<\om$, $z\in\ans{0,1}$, 
and\/ $p_2$ \dd\dP forces 
$\bint\sg\dob(k)=z$ then so does\/ $p_3$. 
\ecl
\bpf[Claim]
Otherwise  
%{\ubf this is not the case}. 
%Then 
there are conditions $p_2$ and $p_3$ as in the claim, 
such that $p_2$ \dd\dP forces $\bint\sg\dob(k)=0$ while 
$p_3$ \dd\dP forces $\bint\sg\dob(k)=1$.
We can wlog assume that
$\dom{u_3}=\dom{u_2}=\text{ some }n$ and  
$m_1\le n\le m$, so overall
\beq
\lam{xx3}
n_0=\dom{u_0}\le m_0 \le n_1=\dom{u_1}\le m_1 \le 
n=\dom{u_2}=\dom{u_3}\le m\,.
\eeq
And we can wlog assume that 
\benq
\nenu
%\setcounter{enumi}{\value{equation}}
%\stepcounter{equation}
\itlb{xx4}\msur
$t_2=t_3 = \text{ some } t\in\zse$, thus  
$p_2=\ang{m,f,u_2,t}$ \dd\dP forces $\bint\sg\dob(k)=0$ while 
$p_3=\ang{m,f,u_3,t}$ (same $m,f,t$!) 
\dd\dP forces $\bint\sg\dob(k)=1$. 
\een
Indeed just let $t=\sup\ans{t_2,t_3}$ termwise, 
thus $t\in\gM$ is a \dd{(\df\cap\gM)}name saying: 
I am a real in $\zn$ and each value $t(j)$ is equal 
to $\sup\ans{t_2(j),t_3(j)}$.

It is clear that the difference between 
the conditions $p_2$ and $p_3$ 
of \ref{xx4} is located in the set 
$U=\ens{j}{u_2(j)\ne u_3(j)}\sq [n_1,n)=\ens{j}{n_1\le j<n}$, 
which we divide into subsets 
$U_2=\ens{j}{u_3(j)<u_2(j)}$ and $U_3=\ens{j}{u_2(j)<u_3(j)}$.
Now define  
$f_2,f_3\in\zn$ as follows:
\beq
\lam{f23}
\left.
\bay{l}
f_3(j)=
\left\{
\bay{lcl}
f(j)+u_2(j)-u_3(j),&\text{whenever} & j\in U_2\\[1ex]
f(j),&\text{otherwise} &
\eay
\right.;
\\[4ex]
f_2(j)=
\left\{
\bay{lcl}
f(j)+u_3(j)-u_2(j),&\text{whenever} & j\in U_3\\[1ex]
f(j),&\text{otherwise} &
\eay
\right.
;
\eay
\right\}
\eeq
so that $f\le f_2$ and $f\le f_3$ termwise, 
the difference between $f,f_2,f_3$ is still located 
in $U\sq[n_1,n)$, and the termwise sums 
${(f_2\res n)}\ar u_2$, ${(f_3\res n)}\ar u_3$ coincide. 

Note that $q_2=\ang{m,f_2}$ and $\ang{m,f_3}$ are 
conditions in $\df\cap\gM$, 
and $f_2\res n_1=f_3\res n_1=f\res n_1$ by construction. 
Let $a_0\in\zn$ be a \dg\ real over $\gM$,
compatible with the condition 
$\ang{m,f}$, so that 
\ben
\aenu
\itlb{yya}
$f\res m\su a_0$ and $f\le a_0$ termwise, 
\een
Accordingly define $a_2,a_3\in\zn$ so that 
\ben
\atc
\aenu
\itlb{yyd}\msur
$a_2\res n=f_2\res n$, 
$a_3\res n=f_3\res n$, 
and $a_3(j)=a_2(j)=a_0(j)$ for all 
$j\ge n$, so that
$f_2\le a_2$ and $f_3\le a_3$ termwise.
\een
Then $a_2,a_3$ are \dg\ reals over $\gM$,
compatible with resp.\  
$\ang{m,f_2}$, $\ang{m,f_3}$.

Now come back to the name $t$ which occurs in conditions 
$p_2,\,p_3$ in \ref{xx4}.
As $t$ is a \dd{(\df\cap\gM)}name for a real in $\zn$, 
in fact the interpretations 
$\bint t{a_0}$, $\bint t{a_2}$, $\bint t{a_3}$ 
belong to $\zn\cap\gM[a_0]$. 
Moreover, as soon as the finite strings 
$f\res n$, $u_2$, $u_3$ (of length $n$) are given, 
the reals $a_2=H_2(a_0)$ and $a_3=H_3(a_0)$ are defined 
by simple functions $H_2$ and $H_3$ whose definitions are 
contained in \ref{yyd} and \eqref{f23}. 
Let $t'\in\gM$ be a \dd{(\df\cap\gM)}name for a real in $\zn$, 
explicitly defined as the termwise supremum of 
$\bint t{\doa}$, $\bint t{H_2(\doa)}$, $\bint t{H_3(\doa)}$, 
so that in particular 
\ben
\atc
\atc
\aenu
\itlb{yyh}\msur
$\bint{t'}{a_0}(j)=
\sup\ans{\bint t{a_0}(j),\bint t{a_2}(j),\bint t{a_3}(j)}$ 
for all $j<\om$. 
\een

Note that 
$q_2=\ang{m,f_2,u_2,t}$ and  
$q_3=\ang{m,f_3,u_3,t}$ are 
still conditions in ${\dP}$, 
and $f_2\res n_1=f_3\res n_1=f\res n_1$ by construction. 
As $n_0\le m_0\le n_1$ by \eqref{xx3}, it follows 
that $p_0\le q_2$ and $p_0\le q_3$. 
(We do not claim that $p_1\le q_{2,3}$ or $p_{2,3}\le q_{2,3}$!)
By the choice of $a_0$ there is a real $b_2\in\zn$ 
such that $\ang{a_0,b_2}$ is a \dd\dP generic pair in 
$\zn\ti\zn$, compatible with the condition 
$p'_2=\ang{m,f,u_2,t'}$, so that 
\ben
\aenu
\atc
\atc
\atc
\itlb{yyb}
$u_2\su b_2$, and 
$u_2\we \bint{t'}{a_0}\le b_2$ 
%, hence $u_2\we t[a_0]\le b_3$, 
termwise.
\een
%Recall that $t[a_0]$ is the \dd ainterpretation of the name $t$.

We further define $b_3\in\zn$ so that 
\ben
\atc
\atc
\atc
\atc
\aenu
\itlb{yyc}\msur
$u_3\su b_3$, and $b_3(j)=b_2(j)$ for all 
$j\ge n=\dom{u_2}=\dom{u_3}$, hence
$u_3\we \bint{t'}{a_0}\le b_3$ termwise by \ref{yyb}.
\een
It follows that $\ang{a_0,b_3}$ is a \dd\dP generic pair, 
compatible with $p_3=\ang{m,f,u_3,t}$.
We conclude by \ref{xx4} that 
\benq
\nenu
%\setcounter{enumi}{\value{equation}}
%\stepcounter{equation}
\itlb{xx6}
$\bint\sg{b_2}(k)=0$ while $\bint\sg{b_3}(k)=1$, thus 
$\bint\sg{b_2}\ne\bint\sg{b_3}$.
\een 

Then the pairs $\ang{a_2,b_2}$ and $\ang{a_3,b_3}$ are 
\dd\dP generic over $\gM$, and we have 
\benq
\nenu
%\setcounter{enumi}{\value{equation}}
%\stepcounter{equation}
\itlb{xx5}
$a_2\ar b_2=a_3\ar b_3$ --- 
therefore $\bint\tau{a_2\ar b_2}=\bint\tau{a_3\ar b_3}$,
\een 
since 
$(a_2\res n)\ar (b_2\res n)=(f_2\res n)\ar u_2=
(f_3\res n)\ar u_3=(a_3\res n)\ar (b_3\res n)$ by construction, 
and if $n\le j$ then $a_3(j)=a_2(j)=a_0(j)$ and $b_3(j)=b_2(j)$. 

Assume for a moment that  
\benq
\nenu
%\setcounter{enumi}{\value{equation}}
%\stepcounter{equation}
\itlb{xx7}
the pairs 
$\ang{a_2,b_2}$, $\ang{a_3,b_3}$ 
are compatible with the conditions resp.\ $q_2,\,q_3$.
\een
%generic pairs 
%$\ang{a_2,b_2}$ and $\ang{a_3,b_3}$ 
%are compatible with the conditions resp.\ $q_2$ and $q_3$. 
Then, as $p_0\le q_2,q_3$, we have 
$\bint\sg{b_2}=\bint\tau{a_2\ar b_2}$ and 
$\bint\sg{b_3}=\bint\tau{a_3\ar b_3}$, 
by \ref{xx1}.
It follows that $\bint\sg{b_2}=\bint\sg{b_3}$ by \ref{xx5}, 
which is 
a contradiction with \ref{xx6}, and this proves the claim.
Thus it remains to establish \ref{xx7}, which 
amounts to
\ben
%\Renu
\def\theenumi{\ref{xx7}\fnsymbol{enumi}}%
\def\labelenumi{\theenumi:}%
\itlb{zza}
$f_2\res m\su a_2$, $f_3\res m\su a_3$,  
and $f_2\le a_2$, $f_3\le a_3$ termwise, 

\itlb{zzb}
$u_2\su b_2$, $u_3\su b_3$, and 

\itlb{zzc}
$u_2\we \bint t{a_2}\le b_2$ and 
$u_3\we \bint t{a_3}\le b_3$ termwise.
\een

Beginning with \ref{zza}, note that $f_2\res n\su a_2$ by 
\ref{yyd}, while if $n\le j<m$ then $a_2(j)=a_0(j)=f(j)$
by \ref{yyd} and \ref{yya}, and $f_2(j)=f(j)$ by construction, 
hence $a_2(j)=f_2(j)$, and $f_2\res m\su a_2$ is verified.
Similarly, if $j\ge m$ then $f_2(j)=f(j)$ and $a_2(j)=a_0(j)$, 
but $f(j)\le a_0(j)$ by \ref{yya}, hence $f_2(j)\le a_2(j)$.

Claim \ref{zzb} immediately follows from \ref{yyb}, \ref{yyc}.

As regards for \ref{zzc}, we have 
$\bint t{a_2}\le \bint{t'}{a_0}$ and 
$\bint t{a_3}\le \bint{t'}{a_0}$ componentwise by \ref{yyh}.
It remains to refer to \ref{yyb} and \ref{yyc}.
\epF{Claim~\ref{xx2}}

A standard consequence of the claim is that $p_1$ 
\dd\dP forces that $\bint\sg{\dob}\in \ugm[\doa]$.
However $p_0\le p_1$ and $p_0$ forces the opposite, 
a contradiction.
\epF{Lemma~\ref{cH2}}

\bpf[Theorem~\ref{mt}, case \ref{mt3}]
As above, 
the {\ubf contrary assumption} 
leads to a formula $\vpi(z)$ with $\ga_0\in\Ord$ as a parameter, 
a condition $p_0=\ang{m_0,f_0}\in\df$ in $\rV$
which \dd\df forces, over $\rV$, that the set 
$\ens{z\in\dn}{\vpi(z)}$ is countable and 
%(by the contrary assumption) 
$\sus z\:(z\nin\uv\land\vpi(z))$,
a sequence $\sis{t_n}{n<\om}\in\rV$ of
\dd\df names for reals in $\zn$, 
and a canonical \dd\df name $T\in\rV$ for 
$\ens{\bint{t_n}{\dot a}}{n<\om}$, 
such that 
\benq
\nenu
\itlb{xx9}
if $x\in\zn$ is a \dg\ real, over $\rV$, 
compatible with ${p_0}$ then 
it is true in $\rV[x]$ that 
$
\ens{z\in\dn}{\vpi(z)}=\ens{\bint{t_n}x}{n<\om}=\bint Tx 
\not\sq\rV\,.
$
\een

Pick a regular cardinal $\ka>\al_0$, sufficiently large for 
$\hk\ka$ to contain $\ga_0$ and all names $t_n$ and $T$.
%and  for Lemma~\ref{aL}\ref{aL0},\ref{aL1} to be true for
%$\gM=\hk\ka$. 
%
Consider a countable elementary submodel $\gM$ of $\hk\ka$ 
containing $\ga_0$, all $t_n$, $T$, and $\df$. 
Let $\pi:\gM\to\gM'$ be the Mostowski collapse onto a 
transitive set $\gM'$.
We have $\pi(t_n)=t_n$ for all $n$ 
(as by the ccc property of $\df$ we can assume that 
$t_n$ is a hereditarily countable set), 
and $\pi(T)=T$. 

By the countability, 
there is a real $a\in\zn$ in $\rV$, \dg\ over $\gM'$.
We can wlog assume that $a(j)=0$ for all $j< m_0$ and 
$a(j)\ge0$ for all $j\ge m_0$.

Let $b\in\zn$ be a real \dg\ over $\rV$, compatible 
with $p_0$. 
In our assumptions, the real $b'=a\ar b\in\zn$  
also is \dg\ over $\rV$ and compatible with $p_0$, 
and $\rV[b']=\rV[b]$ (since $a\in\rV$).
Then $\bint Tb=\bint T{b'}$ by \ref{xx9}. 

On the other hand, $b$ is \dg\ over $\gM'[a]$ as well 
by Lemma~\ref{cI}.
It follows by Lemma~\ref{cH2} that 
$\gM[b]\cap\gM[b']\cap \dn\sq\gM$, therefore
$$
\bint{T}{b}\cap\bint{T}{b'}\sq\gM'[b]\cap\gM'[b']\sq\gM'\sq\rV\,,
$$ 
so that $\bint Tb=\bint T{b'}\sq\rV$, 
and we get a contratiction required with \ref{xx9}.\vom 

\epF{Theorem~\ref{mt}, case \ref{mt3}}

\parf{Sacks case}
\las{sax}

%Here we prove Case \ref{mt4} of Theorem~\ref{mt}.
It is a known property of Sacks-generic extensions 
$\rV[a]$ that if $b\in\dn$ is a real in $\rV[a]$ 
then either $b\in\rV$ or $b$ itself is Sacks-generic 
over $\rV$ and $\rV[b]=\rV[a]$. 
Thus if $X\in\rV[a]$ is an \od\ set of reals in $\rV[a]$ 
and $X\not\sq\rV$ then there is a perfect set $Y\sq\dn$ 
coded in $\rV$, such that every Sacks-generic real 
$b\in Y$ in $\rV[a]$ belongs to $X$. 
However it is true in $\rV[a]$ that every (non-empty) 
perfect set coded in $\rV$ contains uncountably many 
reals Sacks-generic over $\rV$. 

This is a rather transparent argument, 
so we can skip details. \vom 

\qeD{Theorem~\ref{mt}, case \ref{mt4}}

\parf{The Solovay model}
\las{solm}

\bdf
\lam{sm?}
The \rit{first Solovay model\/} is a model of $\zfc$
defined as a generic extension
$\rL[G]$ of $\rL$ by the Levy collapse below an inaccessible
cardinal in $\rL$.
The \rit{second Solovay model\/} is a model of
$\zf+\dc$ equal to the collection of all
hereditarily real-ordinal definable (\hrod)
%, for brevity)
sets in the first model, $\rL[G]$.
\edf

Thus we explicitly consider the case when the ground
\zfc\ model of the Solovay models considered is the
constructible model.
Theorem~\ref{mts} is true for an arbitrary ground model
(with a strongly inaccessible cardinal), but we stick
to the particular case to avoid some minor unrelated
complications.
We'll make use of the following result, implicit
in Stern \cite[proof of 3.2]{stern} and \cite{ksol}.

\bpro
\lam{silt}
It holds in either of the Solovay models, that if an \od\ 
equivalence relation on $\bn$ has at most countably many 
equivalence classes then all of them are 
\od\ sets.\qed 
%, \cite{ksol,stern}. 
\epro

Our first proof of Theorem~\ref{mts}(i) was presented 
in \cite{arXiv1}.
Further research demonstrated though that the proof was 
a largely unnecessary roundabout, and the result can be 
obtained by a rather brief reduction to \ref{silt}. 
We also note that the case,
when $\cX$ is a (non-empty \od\ countable) 
{\bf set of reals} in Theorem~\ref{mts}(i), is well 
known and is implicitly contained in the proof of
the perfect set property for \rod\ sets of reals
by Solovay~\cite{solo}.
Hovever the proofs known for this particular case  
(as, \eg, in \cite{ksol} or Stern~\cite{stern}) 
do not work for sets $\cX\sq\pws{\dn}$.

\bpf[Theorem~\ref{mts}(i)]
\rit{Arguing in the first Solovay model},
let $\cX$ be a 
non-empty \od\ countable {\bf set of sets of reals}; 
we have to prove that $\cX$ 
contains an \od\ element (an \od\ set of reals).
Consider a particular case first.

{\it Case 1\/}: 
$\cX$ consists of \rit{pairwise disjoint} sets of reals. 
If $x,y$ are reals then define $x\rE y$ iff either both $x,y$ 
do not belong to $\bigcup\cX$ or $x,y$ belong to the same 
set $X\in \cX$. 
This is an \od\ equivalence
relation with countably many 
equivalence classes, and hence each \dd\rE class is an \od\ 
set by \ref{silt}, as required. 

{\it Case 2\/}: general. 
Let $\cC$ be the set of all countable sets $C$ of reals, such 
that if $X\ne Y$ belong to $\cX$ 
then already $X\cap C\ne Y\cap C$.
Note that $\cC\ne\pu$ as $\cX$ is countable.
If $X\in\cX$ then let $P_X$ be the set of all pairs of the 
form $\ang{C,X\cap C}$, where $C\in\cC$. 
Then $P_X\cap P_Y=\pu$ whenever $X\ne Y$ belong to $\cX$. 
We conclude that 
$\cP=\ens{P_X}{X\in\cX}$ is a countable collection 
of pairwise disjoint non-empty sets $P_X$ of pairs 
of the form  $\ang{C,C'}$, 
where $C'\sq C$ are countable sets of reals. 

There exists an \od\ coding of such pairs by reals, 
that is, an \od\ map $x\mto \ang{C_x,C'_x}$, where 
$x\in\bn$ is a real, 
$C'_x\sq C_x$ are countable sets of reals for any $x$, 
and for any such pair $\ang{C,C'}$ there is at least 
one $x\in\bn$ such that $C=C_x$ and $C'=C'_x$. 
It follows from the above that the derived sets 
$$
Q_X=\ens{x\in\bn}{\ang{C_x,C'_x}\in P_X}\,,\quad
X\in\cX\,,
$$
form a countable \od\ family $\cQ=\ens{Q_X}{X\in\cX}$
of pairwise disjoint non-empty sets of reals. 
By the result in Case 1, all sets 
$Q_X\in\cQ$ are \od.
But if any $Q_X$ is \od\ then so is 
both $P_X=\ens{\ang{C_x,C'_x}}{x\in Q_X}$ and $X$ itself.

\epF{Theorem~\ref{mts}(i)}

\bpf[Theorem~\ref{mts}(ii)]
\rit{Arguing in the second Solovay model}, let $X\ne\pu$ be
an $\od$ set.
Let $x_0\in X$.
We make use of the fact that, in this model, every set is
real-ordinal definable (\rod).
Thus $x_0$ is \rod; there is an \dd\in formula
$\vpi(\cdot,\cdot,\cdot)$, an ordinal $\al_0$, and a real
$r_0\in\dn$ such that $x_0=F(\al_0,r_0)$, where
$$
F(\al,r)=
\left\{
\bay{rcl}
\text{the only $x$ satisfying $\vpi(\al,r,x)$},&
\text{whenever}& \sus!\,x\,\vpi(\al,r,x)\\[1ex]
\pu,&\text{otherwise}&.
\eay
\right.
$$
Let $R_0=\ens{r\in\dn}{F(\al_0,r)\in X}$, and if
$r,q\in R_0$ then define $r\rE q$ iff
$F(\al_0,r)=F(\al_0,q)$.
Then $\rE$ is an \od\ equivalence relation on an
\od\ set $R_0$.
Moreover $\rE$ has countably many classes
(since $X$ is countable).
It remains to refer to Proposition~\ref{silt}.

\epF{Theorem~\ref{mts}(ii)}

\parf{Problems}
\las{prob}

\bvo
\lam{vo1}
Is the stronger result as in Theorem~\ref{mts}(i) 
(for a set of sets of reals)  
still true in the generic extensions 
mentioned in Theorem~\ref{mt}? 
\evo

\bvo
\lam{vo2}
Is it still true in the first Solovay model that  
every nonempty countable \od\ set (of any kind) 
contains an \od\ element?
\evo

\bvo
\lam{vo3}
Do some other simple generic extensions by a real 
(other than Cohen-generic, Solovay-random, dominating,
ans Sacks) 
admit results similar to Theorem~\ref{mt} 
and also those similar to the old folklore lemmas 
\ref{addc} and \ref{addr} above? 
It would also be interesting to investigate the state of affairs 
in different `coding by a real' models as those defined 
in \cite{kb,klkl}.
\evo

\small
%%%%%%%%%%%%%%

%\bibliographystyle{plain}{\small\bibliography{aml}}

\end{document}